\documentclass[12pt, twoside]{article}
\usepackage{amsmath,amsthm,amssymb}
\usepackage{times}
\usepackage{enumerate}

\usepackage{epsfig}

\def\eqsp{\noalign{\vskip 5pt}}
\def\dsp{\displaystyle}

\newcommand{\lp}{\left (}
\newcommand{\rp}{\right )}
\newcommand{\eqn}[1]{(\ref{#1})}

\theoremstyle{definition}

\numberwithin{equation}{section}

\newcommand{\subjclass}[1]{\bigskip\noindent\emph{2010 Mathematics Subject Classification:}\enspace#1}
\newcommand{\keywords}[1]{\noindent\emph{Keywords:}\enspace#1}

\begin{document}


\baselineskip=17pt


\title{Fine numerical analysis of the crack-tip position for a Mumford-Shah minimizer}

\author{Zhilin Li\\
Center for Research in Scientific Computation \\
and Department of Mathematics\\
North Carolina State University\\
Raleigh, NC 27695, USA\\
zhilin@ncsu.edu\\
Hayk Mikayelyan\\
Mathematical Sciences\\
The University of Nottingham Ningbo\\
199 Taikang East Road, Ningbo 315100, PR China\\
Hayk.Mikayelyan@nottingham.edu.cn}

\date{}

\maketitle


\begin{abstract}

A new algorithm to determine the position of the crack (discontinuity set) of certain
minimizers of Mumford-Shah functional in situations when a crack-tip occurs is introduced.
The conformal mapping $\tilde{z}=\sqrt{z}$ in the complex plane is used to transform the free discontinuity problem to a new type of free boundary problem, where the symmetry of the free boundary is an additional constraint of a non-local nature. Instead of traditional Jacobi or Newton iterative methods, we propose a simple iteration method which does not need the Jacobian but is way fast than the Jacobi iteration. In each iteration, a Laplace equation needs to be solved on an irregular domain with a Dirichlet boundary condition on the fixed part of the boundary; and a Neumann type boundary condition along the free boundary. The augmented immersed interface method is employed to solve the potential problem. The numerical results agree with the analytic analysis and provide insight into some open questions in free discontinuity problems.


\end{abstract}






\section{Introduction}

There are several applications of Mumford-Shah functional introduced in \cite{MS}, but the one which is more relevant to our study is the fracture mechanics. This theory is of course much older and still has many unanswered questions.  It studies complex physical phenomena such as elasticity, plasticity, stress and strain, friction and non-penetration of crack-faces, etc. The most challenging problem though seems to be the simulation and analysis of the situation at the crack-front or, in two dimensional models, crack-tip. This is where one expects the crack to grow further and this is where one wants to predict crack's behavior. The major breakthrough in this field was achieved by an English aeronautical engineer Alan Arnold Griffith in \cite{G}, who was the first to describe rigorously the stress of the fracture and to explain certain effects considered as contradicting the theory by his time.

In 1998 Francfort and Marigo \cite{fm} proposed a quasi-static model based on variational approach. In 2002 Dal Maso and Toader \cite{DMT} obtained existence results based on time-discretization, i.e. for small time steps
to obtain a sequence of crack configurations $(u_n,\Gamma_n)$ by minimizing in the following version of the Mumford-Shah functional
\begin{equation}\label{crack}
J(v,\Gamma) = \int_{B_1\setminus\Gamma} |\nabla v|^2 d x + \lambda^2 \frac{\pi}{2} \, \mathcal{H}^1(\Gamma),
\end{equation}
where $v$ satisfies the Laplace equation on ${B_1\setminus\Gamma}$ and the  given boundary data $v\big|_{\partial B_1}=u_{_D,n}$ at time $t_n$ under additional condition
\begin{equation}\label{addcond}
\Gamma_{n-1}\subset\Gamma
\end{equation}
of the new crack containing the old crack. Here for simplicity we consider the problem in a unit ball $B_1$ and in the formula (\ref{crack}), and
$ \mathcal{H}^1$ denotes the one-dimensional Hausdorff measure (see \cite{AFP}), which for Lipschitz curves coincides with the length.
In order to obtain regularity results for this model one would need better regularity results for
the minimizers of the Mumford-Shah
functional (without additional condition (\ref{addcond})) as those known at this moment, especially results at the crack-tip. From now on we will discuss the minimizers of (\ref{crack}) {\bf without} additional condition (\ref{addcond}) and with given boundary data $u=u_D$.
The authors would like to refer to the excellent monographs by G. David \cite{D} and by L. Ambrosio, N. Fusco and D. Pallara \cite{AFP} for vast amount of knowledge about this well studied problem, which still contains a lot of open questions.

John Andersson and the second author obtained some results about the asymptotics of the Mumford-Shah minimizer
at the crack-tip (see \cite{AM}).
The aim of this paper is the numerical analysis of the problem and precise qualitative and quantitative
investigation of the behavior of minimizers of the Mumford-Shah functional near
the crack-tip. In particular the authors verify numerically the conjecture stated in \cite{AM}, that the curvature at the crack-tip vanishes (see Section~\ref{sec-asymp}).

Different aspects of the regularity theory of Mumford-Shah minimizers have been addressed in several recent publications, such as \cite{DF}, \cite{BM} and \cite{L}.

It is challenging to develop high order method for open-ended interfaces like in the example of crack problems.
In particular it is very hard to work near the crack-tip/crack-front, where in the bulk-term a high order singularity
occurs. In fact the crack-tip is the only place where the two terms of the functional (\ref{crack}), the Dirichlet energy and the length of the crack, scale of same
order. In all other points the crack-length term is dominant, which enables the proof of the known regularity results.

The paper is organized as follows: in Section~\ref{sec-mathform} we present the known results and give the mathematical formulation of the problem, in Section~\ref{sec-nummeth} we describe the numerics and in Section~\ref{sec-numexp} we present and discuss the results.


\section{Mathematical formulation}\label{sec-mathform}

Throughout the paper we are going to use the following function in different coordinate representations (complex, Cartesian, polar)
\begin{equation}\label{sqrt-z}
\Im \sqrt{z}=\frac{\text{sgn}(y)}{\sqrt{2}}\sqrt{\sqrt{x^2+y^2}-x}=r^{\frac{1}{2}}\sin\phi/2,
\end{equation}
where $z=x+iy$, $x=r\cos \phi$, $y=r\sin\phi$ and the discontinuity in (\ref{sqrt-z}) is taken over the branch-cut $(-\infty,0)$ or, depending on the context, over
the set $\Gamma$. In the latter case, the formula in Cartesian coordinates will require slight modification, see (\ref{aux-form}) for
the half-line $\Gamma$ different from $(-\infty,0]$.

It is proven in \cite{BD} that the pair
\begin{equation}
\left (\,\,\lambda r^{1/2}\sin\frac{\phi}{2}\,\,\,;\,\,\, \{(x,0)|-\infty< x\leq 0\}\,\, \right)
\end{equation}
is a global minimizer of the Mumford-Shah functional (\ref{crack}) in the plane ($\lambda$ is the parameter from~(\ref{crack})).
This is understood as being the absolute minimizer in any bounded sub-domain, under its own boundary conditions.
In particular, in the unit disc $B_1$
for the boundary value function
$$
\lambda \sin\frac{\phi}{2},
$$
with discontinuity in the point $(-1,0)$, the pair
\begin{equation}
\left( \,\,\lambda r^{1/2}\sin\frac{\phi}{2}\,\, ;\,\, \Gamma_0\,\, \right), \,\,\,\,\text{with}\,\,\,\,\Gamma_0=\left\{(x,0)|-1\leq x\leq 0\right\},
\end{equation}
is the absolute minimizer of the functional (\ref{crack}), under its own boundary conditions.

In the paper we will discuss the minimizers of (\ref{crack}) in $B_1$ with slightly perturbed boundary data
\begin{equation}\label{bdry-data}
u_D(x,y)=\lambda  \sin\frac{\phi}{2} + \epsilon\chi_{(-\frac{\pi}{2},\frac{\pi}{2})}(\phi) \cos\phi=
 \lambda\frac{ \text{sgn}(y)}{\sqrt{2}}\sqrt{\sqrt{x^2+y^2}-x}+\epsilon\,x^+.
\end{equation}
We will assume the stability of the minimization problem under
those perturbations, i.e. that for small perturbations of $u_D$ the
minimizing crack $\Gamma$ is a curve connecting the
point $(-1,0)\in \partial B_1\cap \mathbb{R}$ with an unknown point $(x_*,y_*)$ inside
the ball $B_{\frac{1}{2}}:=B_{\frac{1}{2}}(0)$ and that this curve is a graph of a function in certain coordinates.

For points $(x_*,y_*)$ in the neighborhood of the origin we will minimize the
functional (\ref{crack}) among cracks which start at $(-1,0)$ and end at $(x_*,y_*)$, by solving relevant Euler-Lagrange equations. Then we analyse
how the position of the crack-tip $(x_*,y_*)$ affects the total energy, the stress intensity  factor and
the asymptotics of the crack near the
crack-tip.

\subsection{Euler-Lagrange conditions}\label{sec-EL}

It is well-known (see \cite{AFP}, \cite{D}) that a minimizer $u$ of  (\ref{crack}) together with the boundary condition
\begin{equation} \label{origDir}
 \left.  u\right |_{\partial B_1} = u_D, \quad \mbox{on} \quad \partial B_1
\end{equation}
satisfies the following four Euler-Lagrange (first order) conditions
\begin{eqnarray}\label{origLap}
  \Delta u  &=& 0  \quad \mbox{in} \quad B_1\setminus\Gamma,  \\  \eqsp
 \eqsp
  \left. \frac{\null}{\null} \partial_\nu u^\pm \right |_{\Gamma} &=& 0  \quad \mbox{on} \quad \Gamma,
  \label{origNeu}
  \\ \eqsp
\lambda^2\frac{\pi}{2}\kappa
& =
 & \left[\left |\nabla u^+ \right|^2 - \left|\nabla u^-\right|^2\right] \Big|_{\Gamma}\quad \mbox{on} \quad \Gamma, \label{orig4eqs}
  \\ \eqsp
|SIF_u(x_*,y_*)| & = & \lambda , \label{origsif}
\end{eqnarray}
where $\nu$ is the normal to crack (discontinuity set) $\Gamma$, $u^\pm$ denotes the values of $u$
on different sides of $\Gamma$, $\kappa$ is the curvature of $\Gamma$ and $(x_*,y_*)$ is the crack-tip.
Further $SIF_u(x_*,y_*)$ is the so-called stress-intensity
factor at the crack-tip. It is known that the solutions of the mixed boundary value problem (\ref{origDir})-(\ref{origNeu})
have the following asymptotics near the crack-tip
$$
u(x,y)=u(x_*,y_*)+  SIF_u(x_*,y_*) r^{\frac{1}{2}}\sin \phi/2+ o(r^{\frac{1}{2}}),
$$
where in this context $(r,\phi)$ are certain radial coordinates with the origin being set in $(x_*,y_*)$.
Thus the condition (\ref{origsif}) says that the coefficient of the $r^{\frac{1}{2}}\sin \phi/2$-term, called
stress intensity factor, cannot be arbitrary.

The conditions (\ref{origLap})
and (\ref{origNeu}) are obtained by the standard variation and conditions (\ref{orig4eqs})
and (\ref{origsif}) by domain variation.

\begin{figure}[htbp]
\begin{minipage}[t]{2.9in}
(a) \\

\epsfysize=1.5in
\centerline{$\qquad$\hbox{\protect\epsffile{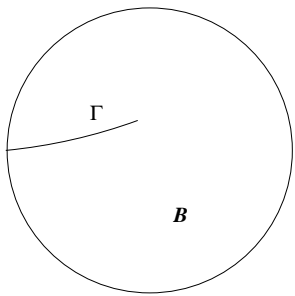}}}
\end{minipage}
\begin{minipage}[t]{2.5in}
(b) \\

\epsfysize=1.6in
\centerline{$\qquad$\hbox{\protect\epsffile{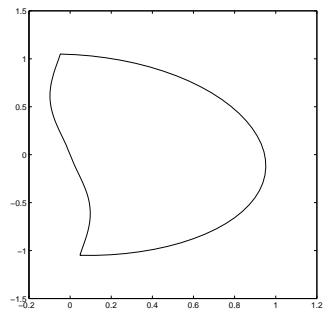}}}
\end{minipage}
\caption{(a): A diagram of a crack $\Gamma$; (b): Domain $\Omega$ after the transformation $\sqrt{z}$.} \label{sqrtz}
\end{figure}

\subsection{Asymptotics at the crack-tip}\label{sec-asymp}

Before stating the result obtained in \cite{AM}, we would like to emphasize that we
use it in our paper only when we numerically verify the conjecture stated in \cite{AM}. The numerical method itself does not rely on this.

\noindent {\bf Theorem.}
{\it
Assume that the minimizer of the functional (\ref{crack}) in $B_1\subset\mathbb{R}^2$
is given by the pair
$(u,\Gamma)$, where $\Gamma=\{(-t,f(t))|0\leq t<1\}$, $f\in C^1([0,1)) \cap C^3((0,1))$,
$f(0)=f'(0)=0$ and $tf''(t)\to_{t\to 0} 0$. Further assume that there exist a limit in $C^1([0,1))$
$$
\lim_{\rho\to 0} \frac{f(\rho t)}{\max_{0<\tau<\rho}|f(\tau)|}\not= 0.
$$
Then there exists a constant $C\not= 0$ such that
\begin{equation}\label{asymptcrack}
f(t)=C t^{1/2+\alpha_k}+o(t^{1/2+\alpha_{k}})
\end{equation}
and
\begin{equation}\label{asymptmt}
u(x,y)=u(0,0)+\lambda \Im \sqrt{z}+ \mathfrak{S}_k(x,y) +C \lambda b_k   r^{\alpha_k} \cos(\alpha_k\phi)+
o(r^{\alpha_k} ),
\end{equation}
where $k<\alpha_k<k+ 1/2$ is one of the positive solutions of
\begin{equation}\label{alf}
\tan (\pi \alpha) =\frac{2}{\pi}\frac{\alpha }{\alpha^2-1/4},
\end{equation}
 $b_k$ is some other absolute constant ($k=1,2,\dots$),
$\mathfrak{S}_1\equiv 0$ and for $k\geq 2$ there are some constants $c_j$ depending on $u$
such that
$$
\mathfrak{S}_k(x,y)=\sum_{j=1}^{k-1} c_j r^\frac{2j+1}{2}\sin\left(\frac{2j+1}{2}\phi\right).
$$
If
$$
\lim_{\rho\to 0} \frac{f(\rho t)}{\max_{0<\tau<\rho}|f(\tau)|}= 0
$$
then $\lim_{t\to 0}t^{-M}f(t)=0$ for any $M>0$.
}

Note that the
first two values of $\alpha_k$ are $\alpha_1 \approx 1.1844$ and $\alpha_2 \approx 2.0989$.
The authors of \cite{AM} have made a conjecture that the coefficient of the $\alpha_1$-term vanishes, and thus in asymptotic expansions
(\ref{asymptcrack}) and  (\ref{asymptmt}) $k\geq 2$. This in particular means that
the curvature of the crack should tend to zero when approaching the crack-tip. In Section~\ref{sec-meas} we present
a numerical justification for this.

\subsection{The $\tilde{z}=\sqrt{z}$ transform}

If we now translate the Figure \ref{sqrtz} (a)
and bring the crack-tip point $(x_*,y_*)$ into the origin and then
use the conformal mapping $\tilde{z}=\sqrt{z}$ with the branch-cut $\Gamma$ in the complex plane, we will obtain a new half-ball like domain $\Omega$ (Figure \ref{sqrtz} (b)), with boundary consisting of two parts; $(\partial \Omega)_D$ originating from $\partial B_1$ and the symmetric $(\partial \Omega)_N$ originating from $\Gamma$.

The new function $\tilde{u}(\tilde{x},\tilde{y})=u(x,y)$ will be defined in $\Omega$ and will satisfy the following conditions equivalent to (\ref{origDir})-(\ref{origsif}),
\begin{equation}\label{Dir}
  \tilde{u} = \tilde{u}_D, \quad \mbox{on} \quad  (\partial \Omega)_D
\end{equation}
where $\tilde{u}_D$ is the boundary data obtained from $u_D$,
\begin{eqnarray}\label{Lap}
  \Delta \tilde{u} &=& 0  \quad \mbox{in} \quad \Omega  \\  \eqsp
   \left. \frac{\null}{\null} \partial_\nu \tilde{u} \right |_{\partial \Omega} &=& 0  \quad \mbox{on} \quad (\partial \Omega)_N
  \label{Neu}  \\ \eqsp
  \pi \lambda^2\left((\tilde{x}^2+\tilde{y}^2)^\frac{1}{2}\tilde{\kappa}(\tilde{x}, \tilde{y})+
\frac{\tilde{x}\tilde{y}'-\tilde{x}'\tilde{y}}{(\tilde{x}^2+\tilde{y}^2)^\frac{1}{2}(\tilde{x}'^2+\tilde{y}'^2)^{\frac{1}{2}}}\right)
& =& \left |\partial_\tau \tilde{u} \right|^2(\tilde{x}, \tilde{y}) - \left|\partial_\tau \tilde{u}\right|^2(-\tilde{x}, -\tilde{y}) \,\,\,\,\,\,\,\,\,\,\,\,\,\,\,\label{4eqs}
  \\ \eqsp
  |\nabla \tilde{u}(0,0)| &= &\lambda,\label{sif}
\end{eqnarray}
where $\nu$ and $\tau$ are normal and tangential directions of the boundary $(\partial \Omega)_N$,
$\tilde{\kappa}$ is the curvature of $(\partial \Omega)_N$.
Observe that the condition (\ref{4eqs}) together with the additional symmetry requirement with respect to the origin imposed on $(\partial \Omega)_N$ have in a sense a non-local nature.

Let us also observe that the asymptotic behavior of the crack near the crack-tip after a rotation will be given by (\ref{asymptcrack}).
This means that after the transformation if the points $(\tilde{f}(t), t)$ describe the upper half of the boundary $(\partial \Omega)_N$, then \begin{equation}\label{asymptopencrack}
\tilde{f}(t)=ct+\frac{C}{2} t^{2\alpha_k}+o(t^{2\alpha_{k}}).
\end{equation}

As we see for a fixed $(x_*,y_*)$ solving the Mumford-Shah minimization problem is equivalent
to solving the over-determined mixed (Dirichlet-Neumann) boundary value problem (\ref{Dir})-(\ref{4eqs})
with and additional symmetry constraint imposed on $(\partial \Omega)_N$.

After solving (\ref{Dir})-(\ref{4eqs}) for each $(x_*,y_*)$ we will compute the Mumford-Shah energy, which will be a function depending on $(x_*,y_*)$. The optimal position of the crack-tip can now be determined by
minimization of that function. We will see that the condition (\ref{sif}) will be satisfied for many values of
$(x_*,y_*)$ lying on an interface, which contains the optimal crack-tip, but the minimizer will be unique. This indicates that the known
Euler-Lagrange conditions are not complete and there is another first order condition missing. Moreover, we will verify numerically that
for the minimizer in the asymptotic expansion (\ref{asymptopencrack}) the value $k>1$, and thus the curvature converges to zero when approaching the crack-tip.


\section{The numerical method}\label{sec-nummeth}

As usual, the key to solve the minimization problem is to use some sort of steepest descent direction plus some acceleration techniques. In each iteration, we need to solve the Laplacian equation with a Dirichlet boundary condition along the fixed part and homogeneous Neumann boundary condition along the free part of the boundary. Note that the Jacobian matrix is difficult and expensive to get.
Our numerical method is based on the augmented immersed interface method \cite{li-ito} for solving Poisson equations on irregular domains so that a fast Fourier transform (FFT) based fast Poisson solver can be utilized. A preconditioning technique is also developed to solve the related Schur complement of the system of equations. To accelerate the convergence of the steepest descent direction method, we use some analysis to obtain a good approximation of the boundary of the minimizer; and a simple but efficient iterative method for the non-linear iteration.
Below we explain our algorithm in detail.

\subsection{Set up the problem}

We start with a rectangular domain $[x_{min}\; x_{max}]\times [y_{min}\; y_{max}]$ that is large enough to possible variations of the free boundary.
The fixed boundary is described as part of the unit circle $r=1$. In discretization, for convenience, we assume that $x_{max}-x_{min}=y_{max}-y_{min}$. Given a parameter
 $n$, we set up a uniform Cartesian grid with $h=(x_{max}-x_{min})/n$,
\begin{equation}\label{grid}
    x_i = x_{min} + i h; \quad y_j = y_{min} + jh, \qquad i,j = 0, 1, \cdots, n.
\end{equation}

Start with an initial tip position $(x_*,y_*)$, we numerically find the function $u$ with the discontinuity set (crack) $\Gamma$, connecting the discontinuity point $(-1,0)$ on the boundary with $(x_*,y_*)$, that minimizes the MS energy.
In order to set the crack-tip in the origin we translate the picture by the vector $(-x_*,-y_*)$ and discretize the fixed boundary with a parameter $N_b$. We set $\Delta s = 2 \pi/N_b$ and a Lagrangian grid
\begin{equation}\label{grid2}
    X_k = -x_* + r \cos \lp (k-1) \Delta s-\pi \rp, \quad Y_k = - y_* +  r \sin \lp (k-1) \Delta s-\pi \rp, \qquad k=1, 2, \cdots N_b,
\end{equation}
where $r=1$ in our simulation. Note that we use $k \Delta s-\pi$ instead of $k \Delta s$ so that the starting and ending point is $(-1-x_*,-y_*)$, the discontinuity point of the boundary data.
In general we should also choose $N_b$ in such a way that $\Delta s \sim h$. 

After the transformation $\tilde{x}+i\tilde{y}=\tilde{z}=\sqrt{z}$ with a branch-cut over the ray starting in the origin and going through $(-1-x_*,-y_*)$ we obtain
\begin{eqnarray}
  \tilde{X}_k  &=& \displaystyle \frac{\text{sgn}(Y_k)\, \text{sgn}_\Gamma(X_k,Y_k)}{\sqrt{2}} \sqrt{\sqrt{X_k^2+Y_k^2}+X_k}, \\ \eqsp
  \tilde{Y}_k  &=&  \frac{\text{sgn}_\Gamma(X_k,Y_k)}{\sqrt{2}} \sqrt{\sqrt{X_k^2+Y_k^2}-X_k},\label{aux-form}
\end{eqnarray}
where the function
$$
\text{sgn}_\Gamma(x,y)=
\begin{cases}
\text{sgn}(y) & \text{if} \,\,\, x\geq 0, \\ \eqsp
1 &\displaystyle \text{if} \,\,\, x< 0 \,\,\,\text{and}\,\,\, y\geq \frac{y_*}{1+x_*}x, \\ \eqsp
-1 & \displaystyle\text{if}\,\,\, x< 0\,\,\,\text{and}\,\,\,   y<\frac{y_*}{1+x_*}x ,
\end{cases}
$$
is needed because we have taken a branch-cut different from negative half-axis.


\subsection{Set up an initial free boundary}\label{sec-setup}

From now on, we will use the transformed domain $\Omega$ unless stated otherwise.
By construction the free boundary should be symmetric with respect to the origin $(0,\,0)$. In our iterative process in updating the free
boundary, we just need to compute the upper half and get the lower half from the symmetry. So now we need to
make a good initial guess for $(\partial \Omega)_N$. We start with a straight line connecting
$(\tilde{X}_{N_b}, \tilde{Y}_{N_b})$ and the origin $(0,\,0)$ with the step size $h_1=d/\mbox{int}(h/d)$, where $d$ is the distance between the two points. The straight line can be written as
\begin{eqnarray}
  \tilde{X}_k  &=& a_1 t_k, \qquad \mbox{with} \quad t_0=0, \\
  \tilde{Y}_k  &=&  a_2 t_k
\end{eqnarray}
where $(a_1, a_2)$ is the direction connecting $(0,\,0)$ and $(\tilde{X}_{N_b}, \tilde{Y}_{N_b})$.

\begin{figure}[htbp]
\begin{minipage}[t]{2.3in}
\epsfysize=2.2in
\centerline{\hbox{\protect\epsffile{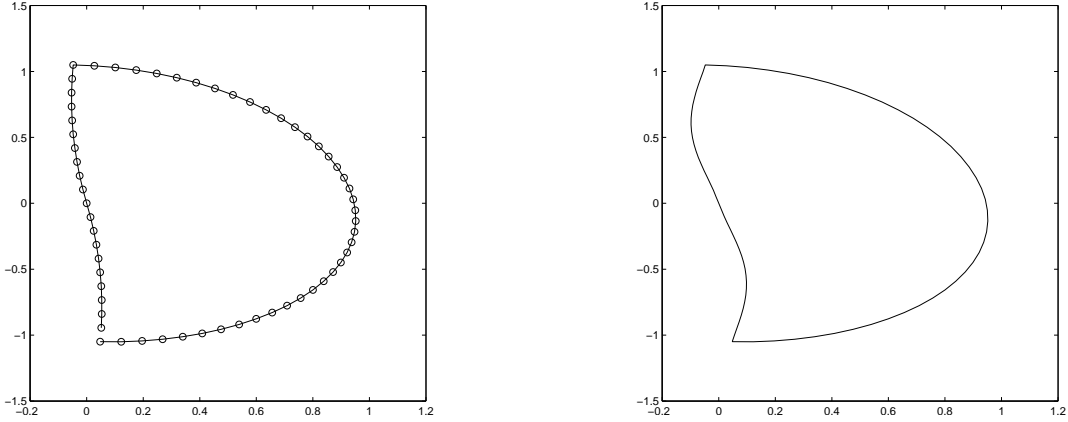}}}
\end{minipage}
\hfil
\begin{minipage}[t]{2.1in}
\epsfysize=2.2in
\centerline{\hbox{\protect\epsffile{plot071014.eps}}}
\end{minipage}
\caption{(a): An initial set-up of the free boundary with one choice $(x_*,y_*)$. The little circles 'o' are discrete points of the boundary called the control points. (b): The final boundary that minimizes the Mumford-Shah energy corresponding to the initial configuration of the left plot with  $\epsilon=10^{-2}$.} \label{fig:point}
\end{figure}

From the theoretical analysis we know that if the angle between the free boundary $(\partial \Omega)_N$ and the fixed boundary $(\partial \Omega)_D$ is not $90$ degrees then in the Dirichlet-Neumann mixed boundary value problem we can have a singularity in the corner. Thus we make a better initial guess by taking a parabola connecting the two points $(\tilde{X}_{N_b}, \tilde{Y}_{N_b})$ and the origin $(0,\,0)$ such that the angle between free and fixed boundaries is 90 degree. We define
\begin{equation}\label{gt0}
    g(t) = c t (t-d), \qquad c = \frac{y_*}{d (1+x_*)}.
\end{equation}
The initial free boundary is then determined by
\begin{eqnarray}
  \tilde{X}_k  &=:& \tilde{X}_k + a_1 g(t_k), \qquad \mbox{with} \quad t_0=0, \\
  \tilde{Y}_k  &=:&  \tilde{Y}_k - a_2 g(t_k),
\end{eqnarray}
here the notation $=:$ stands for overwriting. The bottom part is determined using the symmetry. We show an initial domain with
$x_*=-0.1$, $y_*=0.1$ in Figure~\ref{fig:point}. Note that the angle of the free boundary and the ``half-circle'' is $\pi/2$.

\subsection{Solve the Laplace equation on an irregular domain}

The bulk of computational cost is to solve the Laplace equation on an irregular domain similar to the one in Fig.~\ref{fig:point}~(a) at each iteration.
We use the augmented immersed interface method (AIIM) \cite{li-ito} to solve the Laplace equation with Dirichlet boundary condition along the fixed boundary and homogeneous Neumann boundary condition along the free boundary. The consideration to use the AIIM are two folds. The first one is that we can utilize a FFT based fast Poisson solver. The second one is our experiences using the IIM to solve PDEs on irregular domains with second order accuracy. The difficulty here is the different boundary conditions on different parts of boundaries.  We use the closed cubic spline package \cite{li-pack} to represent the boundary, see Fig.~\ref{fig:point}~(a) for an illustration. First of all, we imbed the domain $\Omega$ into the rectangular box $R=[x_{min}\; x_{max}]\times [y_{min}\; y_{max}]$. The Laplace equation is extended to the entire domain to form an augmented interface problem
\begin{eqnarray}
  && \Delta U = 0 \quad (x,y) \in R, \qquad U_R=0;  \\ \eqsp
  &&  \left [U \right ]_{(\partial \Omega)_D} = 0; \qquad  \left[ \frac{\partial U}{\partial \nu} \right ]_{(\partial \Omega)_D}=q_D, \qquad (X,Y) \in (\partial \Omega)_D, \\ \eqsp
  &&  \left [U \right ]_{(\partial \Omega)_N}  = 0; \qquad  \left[ \frac{\partial U}{\partial \nu} \right ]_{(\partial \Omega)_N}  = q_N, \qquad (X,Y) \in (\partial \Omega)_N, \\ \eqsp
  && \left. U \right |_{(\partial \Omega)_D}  = U_D, \quad \mbox{$U_D$ is known}, \qquad \left.  \frac{\partial U}{\partial \nu} \right |_{(\partial \Omega)_N}  = 0, \label{int_bc}
\end{eqnarray}
where $\nu$ is the unit normal direction of the domain $\Omega$, $[v]$ stands for the jump of a quantity of $v$ across the boundary, thus $[U]_{D}$ means the jump of the solution $U$ across the boundary where the Dirichlet boundary condition is defined. In the system of above, the Dirichlet boundary condition $U_D$, the domain $\Omega$, and the free boundary $(\partial \Omega)_N$, are known; while $U$, $q_N$, and $q_D$ are unknowns. In the AIIM, the augmented variables $q_N$ and $q_D$ should be chosen such that the boundary conditions $U |_{D} = U_D$ and $\frac{\partial U}{\partial n} |_{N} = 0$  are satisfied along with the Laplace equation on $\Omega$. The value of $U$ in $R\setminus \Omega$ is of no interest. We set homogeneous Dirichlet boundary condition along the boundary of the rectangle. Note that to our best knowledge, this is the first time that the augmented IIM has been applied to Laplace equations on irregular domain with different boundary conditions on different parts of the boundary. The boundary conditions, the augmented variables $q_N$ and $q_D$, are defined at the discrete points as explain above.

We use the standard centered five point finite difference stencil to discretize the Laplace equation with IIM correction terms at the right hand side given the discrete values of the augmented variable $q_N$ and $q_D$. The discreization in the matrix-vector form can be written as
\begin{equation}\label{aug1}
    A U + B Q = F_1,
\end{equation}
where $A$ is the matrix corresponding to the discrete 5-point Laplacian, $U$ is the approximate solution defined at all grid points $(x_i,y_j)$, $Q$ is the vector formed from the discrete augmented variable $q_N$ and $q_D$ defined at all boundary points $(X_k, Y_k)$. The solution $U$ should also satisfy the internal boundary condition \eqn{int_bc} at the internal boundary points $(X_k, Y_k)$. In the AIIM,  least squares interpolation schemes are used to approximate the internal boundary condition \eqn{int_bc} at $(X_k, Y_k)$. In the matrix-vector equation, the discretization can be written as
\begin{equation}\label{aug2}
    C U + D Q = F_2.
\end{equation}
In equation \eqn{aug1}-\eqn{aug2}, the matrix $A$ is fixed, while the sparse matrix $B$, $C$, and $D$ are changing with the free boundary. The Schur complement for $Q$ is
\begin{equation}\label{schur}
    \left ( D- C A^{-1} B \right ) Q= F_2-CA^{-1}F_1, \quad \mbox{or} \quad \mathcal{S} \, Q= \mathcal{F}.
\end{equation}
It is time and memory consuming to form those matrices and the Schur complement at each iteration. Instead, we use the GMRES iterative method to solve for $Q$ so that a fast Poisson solver can be utilized. Those sparse matrices $A$, $B$, $C$, and $D$ correspond to finite difference discretization, and interpolation schemes and do not formed explicitly. We refer the readers to \cite{li-ito} for the detailed implementation. The GMRES iterative method only requires the matrix-vector multiplication which contains the two steps: (1) solve for $U$ from $AU = F_1 - B Q$; (2) compute the residual of the boundary condition $R(Q) = C U + D Q - F_2$. We skip the details here since they can be found in \cite{li-ito}.

\subsection{A preconditioning strategy}

Since both Dirichlet and Neumann boundary conditions are prescribed, in general, we do not have the fast convergence of the GMRES method. The number of GMRES iterations depends on the geometry and the mesh size. In \cite{xia-li}, a sophisticated precondition technique is developed for the augmented IIM. In this paper we use a simplified version of the preconditioning strategy based on the idea from \cite{xia-li}.

Let the dimension of the Schur complement matrix be $N_3$.
In the simple preconditioning strategy, we take a parameter $L$,  an integer between $10\sim 30$ such that $N_3/L$ is close to an integer. Let $I_L$be the matrix of $L$ rows and $N_3$ columns whose entries are number one. Define $P=\mathcal{S} I_L$. The preconditioning matrix is defined as
\begin{equation}\label{pre}
    P_c = \left [ P_1^{-1},  P_2^{-1}, \cdots \right ]^T,
\end{equation}
where $P_1^{-1}$ is the inverse of the matrix formed by the first $L$ rows, $P_2^{-1}$ is the inverse of the matrix formed by the next $L$ rows, and so on. The new system of equations is $P_c \mathcal{S} \, Q= P_c\mathcal{F}$. In Table~\ref{tab:pre}, we list a comparison of the number of GMRES iterations from one of examples. We can see the number of GMRES iterations is significantly reduced with suitable choice of $L$.

\subsection{A new iterative scheme for the free boundary problem}\label{sec-iterative}

In order to find the correct crack $\Gamma$ for the fixed crack-tip $(x_*,y_*)$ we need to solve the non-linear boundary value problem (\ref{Dir})-(\ref{4eqs}). For a given boundary  $(\partial \Omega)_N$ the problem (\ref{Dir})-(\ref{Neu}) is well-posed and
can be solved as a mixed type boundary value problem. But there is no guarantee that the condition (\ref{4eqs})
will be satisfied. We will find the correct $(\partial \Omega)_N$ satisfying (\ref{4eqs}) iteratively.

After rotation the part of the curve $(\partial \Omega)_N$ in the upper half space can be modeled by $(t,g(t))$, where $t\in(0,d)$, $d=((1+x_*)^2+y_*^2)^{1/4}$ is the distance between the origin and the contact point of $(\partial \Omega)_D$ and $(\partial \Omega)_N$, and $g(0)=g(d)=0$. Note that here the power $\frac{1}{4}$ is a result of the $\sqrt{z}$ transform. Thus the condition (\ref{4eqs}) can be written as
\begin{equation}\label{fbeq}
\frac{(t^2+g^2(t))g''(t)}{(1+g'^2(t))^\frac{3}{2}} + \frac{tg'(t)-g(t)}{ (t^2+g^2)^\frac{1}{2}(1+g'^2(t))^{\frac{1}{2}}  }
=G(t),
\end{equation}where $t\in (0,d)$ and
$$
G(t)=\frac{1}{\pi\lambda^2}\left(\left |\partial_\tau \tilde{u} \right|^2(t, g(t)) - \left|\partial_\tau \tilde{u}\right|^2(-t, -g(t)).
\right),
$$
with boundary conditions 
$$
g(0)=g(d)=0.
$$
We propose the following iterative scheme to solve (\ref{fbeq}):
\begin{equation}
\label{iterfbeq}
\frac{(t^2+g_n^2(t))}{(1+g'^2_n(t))^\frac{3}{2}}g_{n+1}''(t) + \frac{tg_n'(t)-g_n(t)}{ (t^2+g_n^2)^\frac{1}{2}(1+g_n'^2(t))^{\frac{1}{2}}  }
=G_n(t),
\end{equation}
where
$$
g_{n+1}(0)=g_{n+1}(d)=0,
$$
$$
G_n(t)=\frac{1}{\pi\lambda^2}\left(\left |\partial_\tau \tilde{u}_n \right|^2(t, g_n(t)) - \left|\partial_\tau \tilde{u}_n\right|^2(-t, -g_n(t))
\right). \label{ite-ode}
$$
and $\tilde{u}_n$ is the solution to the mixed boundary value problem (\ref{Dir})-(\ref{Neu}) in the domain
with boundary determined by $(t,g_n(t))$.

As already indicated in Section~\ref{sec-setup} we start the iteration with initial guess $g_1(t)=g(t)$ as in~(\ref{gt0}).

\subsection{A summary of the  algorithm}

The proposed algorithm is to find the variation of $g(t)$ in \eqn{iterfbeq} such that \eqn{fbeq} is satisfied. The algorithm can be outlined below.
\begin{itemize}
  \item Step 1: Set-up the problem. Define a rectangular domain $[x_{min}\; x_{max}]\times [y_{min}\; y_{max}]$. Define $(x_*,y_*)$, $\epsilon$, and other parameters.
  \item Step 2: Determine an initial $g_0(t_k)$ as explained in Section~\ref{sec-setup}. For $l=1, \cdots$ until converges,
  \begin{enumerate}
    \item Solve the Laplace equation with the free boundary fixed.
    \item Solve the ODE \eqn{ite-ode}.
    \item Update the new $g_{l+1}$ and compare $\max_k\| g_{l+1}(t_k) - g_l(t_k)\|$ for the convergence.
  \end{enumerate}
  \item Step 3: Data and visualization analysis.
\end{itemize}

\section{Numerical experiments}\label{sec-numexp}

There are two purposes for the simulations. The first one is to compare to the theoretical analysis. The second purpose is to test the capability of the algorithm to find the crack location that minimizes the Mumford-Shah energy.

\subsection{Accuracy check of the Laplacian solver on irregular domains with both Dirichlet and Neumann boundary conditions}

We first want to make sure that the Poisson solver for the Laplace equation with an irregular domain and  with both Dirichlet and Neumann boundary conditions works properly and accurately.  We use an example for which we know the exact solution to check the method. The boundary is the half circle $x^2 + y^2 \le 1$ and $x\ge 0$ which is imbedded in a large rectangle $R=[-2,\; 2]\times [-2,\; 2]$. The analytic solution is
\begin{equation}\label{poisson-ex}
  u(x,y) = e^{-y}\, \cos x
\end{equation}
which satisfies Laplace equation and the homogeneous Neumann boundary condition along the $y$-axis. The Dirichlet boundary condition is applied according to the exact solution. In Table~\ref{tab-poisson}, we list the result of a grid refinement analysis. We can see that the method is second order accurate in the infinity norm. In the table, the error is defined as
\begin{equation}
\| E\|_{L^\infty}=\dsp \max_{x_{ij}\in\Omega} \left | u(x_i,y_j)- U_{ij} \right|,
\end{equation}
where $U_{ij}$ is the computed solution at the grid point $(x_i,y_j)$ inside the domain $\Omega$. The order of convergence is measured as usual
\begin{equation}\label{conv}
  r = \log \left(E_N/\log E_{2N} \right ) /\log 2,
\end{equation}
where $N$ is the number of grid lines in each coordinate direction. Thus the mesh size is $h_x=h_y = 4/N$. The number of unknowns of $U_{ij}$ is $O(N^2)$ while  the number of augmented variables $Q_i$ is $O(N)$.

\begin{table}[hptb]
\begin{center}
\begin{tabular}{|c|c|c| } \hline
  $N$ &  $E_N$ &  $r$ \\ \hline
  $40$ &  $7.8049e-002 $ &  $ $ \\ \hline
  $80$ &  $2.9358e-002 $ &  $  1.4106   $ \\ \hline
  $160$ &  $3.9157e-003 $ &  $ 2.9064$ \\ \hline
  $320$ &  $1.5625e-004$ &  $ 4.6473$ \\ \hline
  $640$ &  $ 3.8137e-005$ &  $ 2.0346 $ \\ \hline
\end{tabular}
\caption{A grid refinement analysis to check the convergence of the Laplace solver. On average second order accuracy in the $L^{\infty}$ norm is observed. \label{tab-poisson}}
\end{center}
\end{table}

In Table~\ref{tab:pre}, we list number of GMRES iterations with a preconditioning parameter $L$ for the case when $N=640$. The condition number for the Schur complement is $5.0619\times 10^{5}$.  We use $L=0$ for the case with no preconditioning. The restart parameter is
$RESTART=320$ and the tolerance is takes as $10^{-7}$. We can see the simple preconditioning technique does reduce the number of GMRES iterations significantly.
\begin{table}[hptb]
\begin{center}
\begin{tabular}{|c|c|c| c|c|c|c|c|} \hline
  $L$   & $0$   & $5$   & $10$ & $15$ & $20$ & $25$ & $30$ \\ \hline
  $No.$ & $1901$ & $152$ & $86$ & $62$ & $55$ & $48$ & $46$ \\ \hline
 \end{tabular}
\caption{The number of GMRES with and without ($L=0$) the preconditioning technique. \label{tab:pre}}
\end{center}
\end{table}

\subsection{Crack simulations, measurements at $(x_*,y_*)$ and observations}\label{sec-meas}

We start with a point $(x_*, y_*)$, and an initial guess of the crack location as described in Section~\ref{sec-iterative}, see also Figure~\ref{crack}.
We run our algorithm until $\|\Gamma^{m+1}-\Gamma^{m}\|\le tol$, where $tol=10^{-6}$. The final $\Gamma^{m}$ is consider the crack location that minimizes the Mumford-Shah energy among all curves with crack-tip at $(x_*,y_*)$.
First we validate our method for parameters $\lambda=1$, $\epsilon=0$
and the boundary data (\ref{bdry-data}). For $(x_*, y_*)=(0,0)$ the final
free boundary is the line segment $[-1, 1]$ on the $y$-axis.
This numerically confirms the result in \cite{BD}.

Then for the boundary data (\ref{bdry-data}) with $\lambda=1$ and $\epsilon=0.01$ we vary
$(x_*, y_*)=(-0.2 + i_1 \tilde{h}, -0.2 + j_1 \tilde{h})$, where $\tilde{h}=0.02$ and
$i_1, \, j_1=0,1, \cdots, 20$, and determine iteratively the correct free boundary.

\begin{figure}[htbp]
\begin{minipage}[t]{2.3in}
\epsfysize=2.1in
\centerline{$\qquad$\hbox{\protect\epsffile{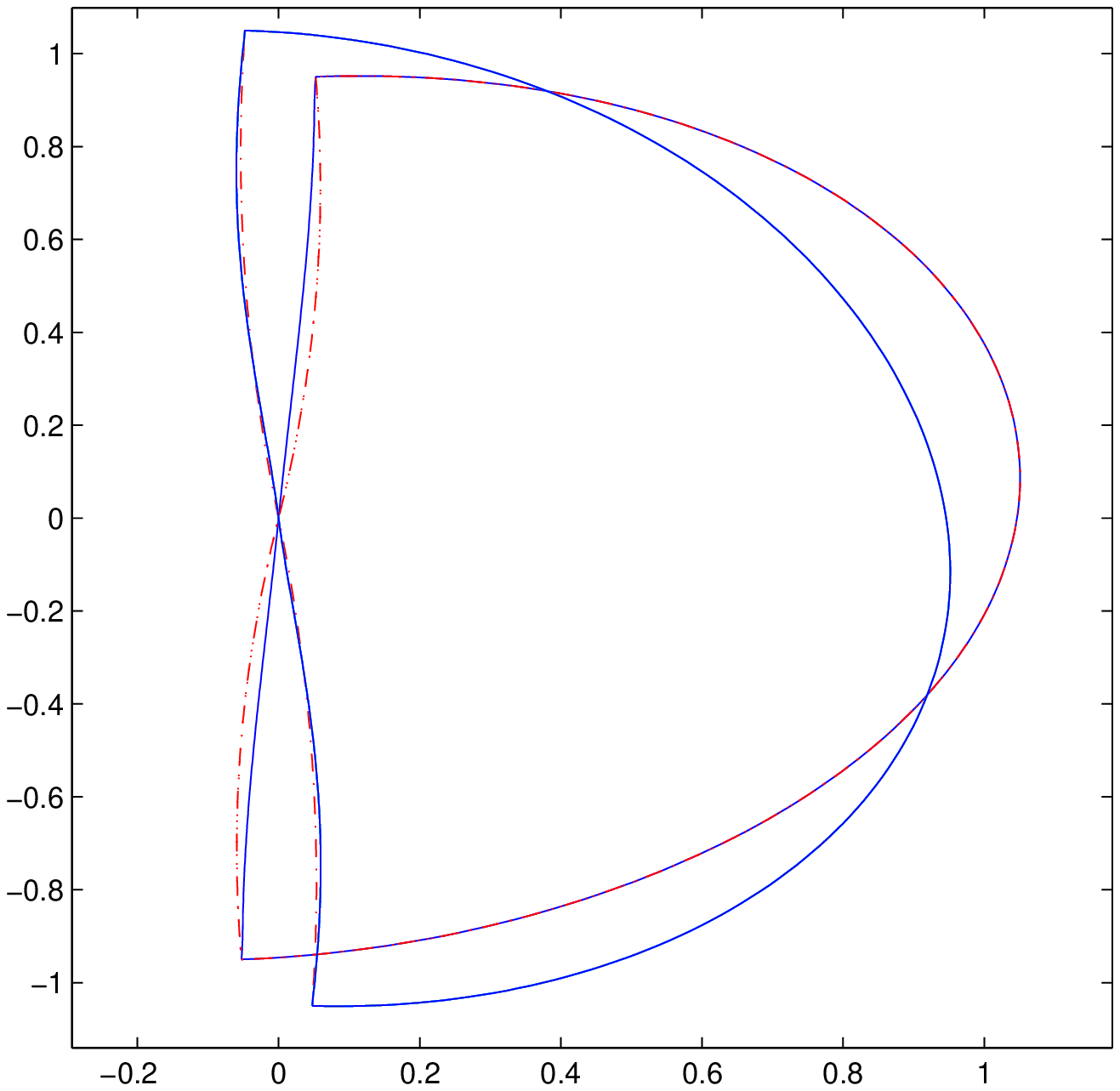}}}
\end{minipage}
\hfil
\begin{minipage}[t]{2.1in}
\epsfysize=2.1in
\centerline{\hbox{\protect\epsffile{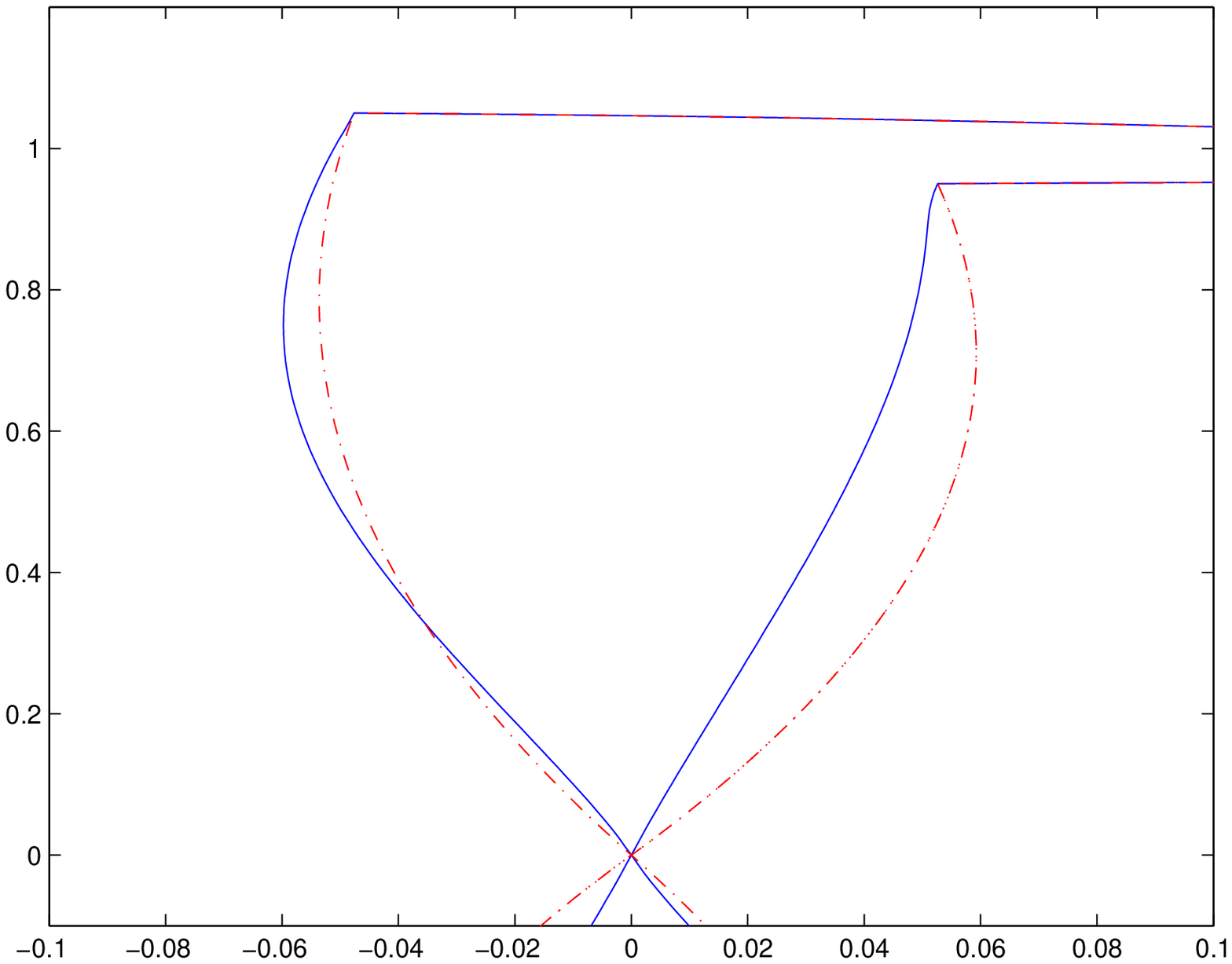}}}
\end{minipage}
\caption{Left plot: the initial (red) and final boundary (blue) of
two different  $(x_*,y_*)$s. One corresponds to $(x_*,y_*)=(0.1,0.1)$;
the other $(x_*,y_*)=(-0.1,-0.1)$. Right plot: zoom-out plot of the
free boundary. Note that the right part of boundary is fixed. } \label{fig:frbdr}
\end{figure}

In Fig.~\ref{fig:frbdr}, we plot the  initial and final boundary of
two different  $(x_*,y_*)$s. One corresponds to $(x_*,y_*)=(0.1,0.1)$;
the other $(x_*,y_*)=(-0.1,-0.1)$. In the right plot of  Fig.~\ref{fig:frbdr},
we zoom out the final free boundary. Note that the right part of the boundary
is the ``half-circle'' $(\partial \Omega)_D$ which is fixed.

As soon as we have the correct solution we compute following quantities as functions depending on $(x_*,y_*)$:

(a)  $MS(x_*,y_*)$ -- the Mumford-Shah energy, which consists of Dirichlet energy of $u$ and the length of the set $\Gamma$. We observe that the Dirichlet energy of the function $\tilde{u}$ equals the Dirichlet energy of the function $u$ in the original picture (before applying $\sqrt{z}$ transform), and can be computed with very high accuracy, since $\nabla \tilde{u}$ has no singularity. The length of the discontinuity set $\Gamma$ connecting the point $(-1,0)$ with $(x_*,y_*)$ can be easily computed after transforming the optimal free boundary $(\partial \Omega)_N$ by the inverse transform $z=\tilde{z}^2$.

\begin{center}
\begin{figure}[htbp]\label{fig:MSE}
\begin{minipage}[t]{6in}
\epsfysize=3.9in
\centerline{$\qquad$\hbox{\protect\epsffile{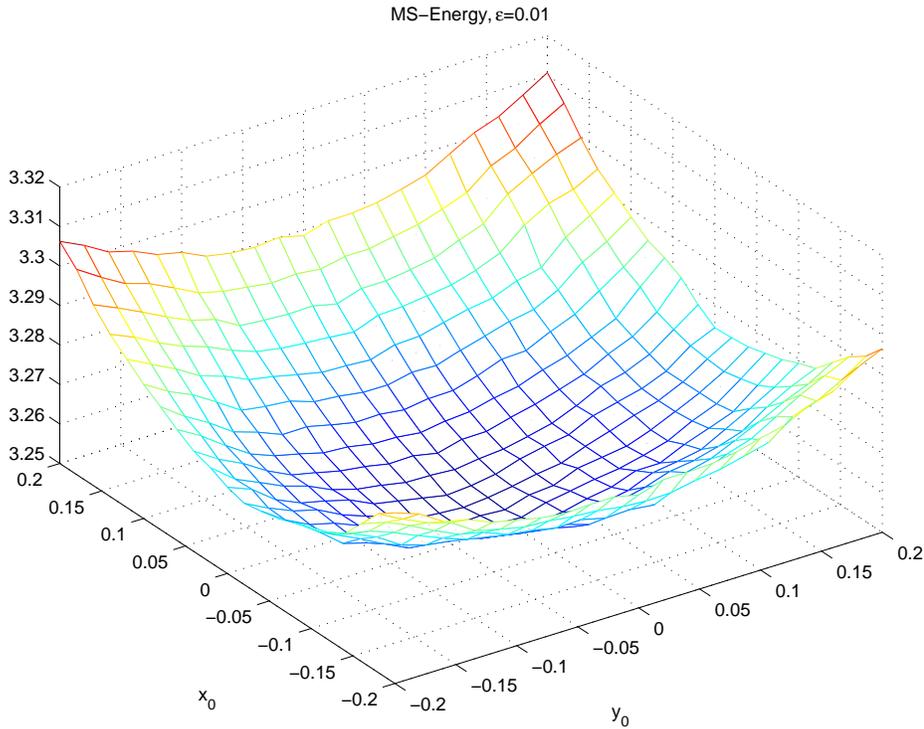}}}
\end{minipage}
\caption{Plot of the MS energy as functions of $(x_*,y_*)$.} \label{fig:ms}
\end{figure}
\end{center}

(b) $SIF_u(x_*,y_*)$ -- the stress intensity factor. As indicated in Section~\ref{sec-EL} $SIF_u(x_*,y_*)=|\nabla \tilde{u}(0,0)|$, which we can measure numerically. Obviously it is easier and more accurate to compute the gradient of a non-degenerate function $\tilde{u}$, rather than the coefficient of the asymptotic term $\Im \sqrt{z}$ of a discontinuous function $u$ with exploding gradient.

(c) the value $u(x_*,y_*)=\tilde{u}(0,0)$.

(d) From (\ref{asymptopencrack}) we know that
$$
g(t)=c_1 t + c\, t^{2\alpha_{k}}.
$$
In order to answer the question whether $k>1$ or not, we apply the non-linear
least square fitting to approximate the solution $g(t)$ by functions of the form
$$
c_1 t + c_2 t^{2\alpha_1}+c_3 t^{2\alpha_2}
$$
in the interval $(0,1/2)$, and obtain the coefficients $c_1, c_2, c_3$ as functions of $(x_*,y_*)$. The
crack-tips $(x_*,y_*)$, where the value of the $c_2$ vanishes should correspond to the configurations with
$k>1$ in the asymptotic expansions (\ref{asymptcrack}) and (\ref{asymptmt}), and thus vanishing curvature at the crack-tip.

\vspace{5mm}

We observe the following results:

(i) The total energy $MS(x_*,y_*)$ is a convex function with minimum at $(-0.02,-0.02)$.

(ii) The function $SIF(x_*,y_*)$ is nearly linear and
the set $SIF(x_*,y_*)=1$
is a vertical interface corresponding to crack-tips with optimal stress intensity factor. This means that there is a continuous family of crack-configurations, with crack-tips on the mentioned interface, which satisfy all known Euler-Lagrange conditions (\ref{origLap})-(\ref{origsif}), but only one those crack-configurations corresponds to an energy minimizer. This is
a numerical observation for existence of another first order condition for Mumford-Shah minimizers.

\begin{figure}[htbp]
\begin{minipage}[t]{2.3in}
  \epsfysize=2.1in
  \centerline{$\qquad$\hbox{\protect\epsffile{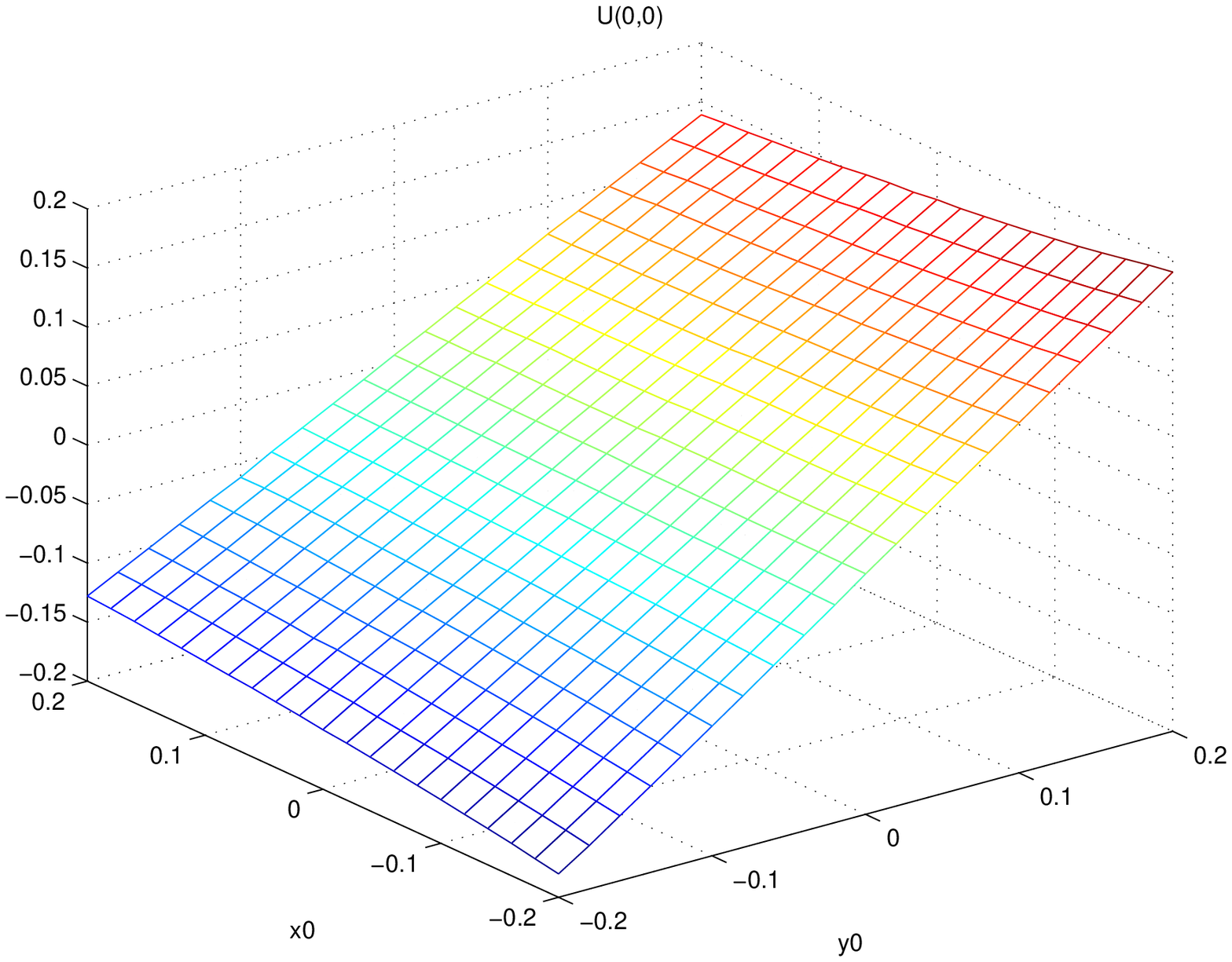}}}
  \end{minipage}
\hfil
\begin{minipage}[t]{2.1in}
  \epsfysize=2.1in
\centerline{$\qquad$\hbox{\protect\epsffile{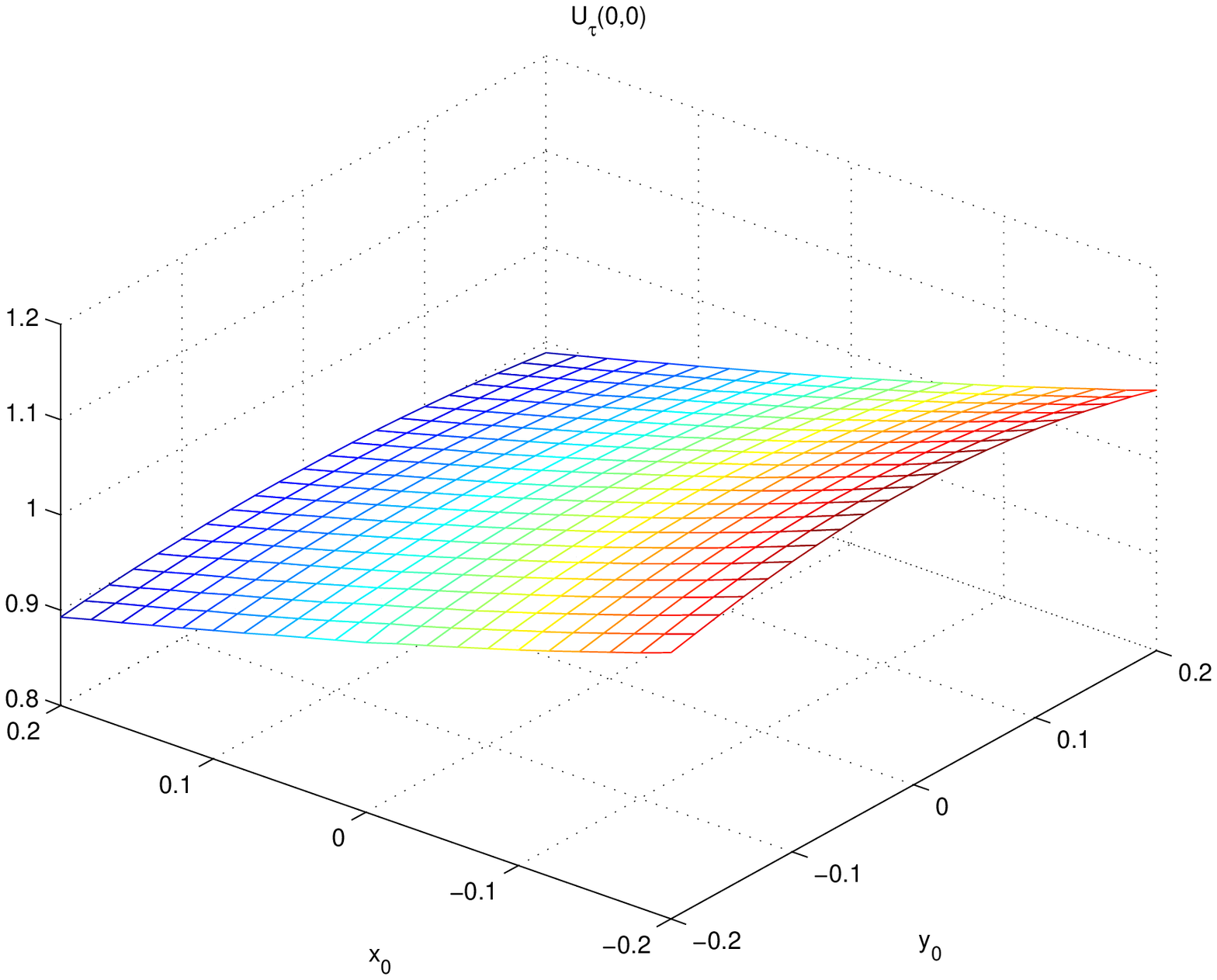}}}
\end{minipage}
\caption{Plot of $u(x_*,y_*)=U(0,0)$ and $SIF(x_*,y_*)=U_{\tau}(0,0)$.}    \label{fig: U00_UT00}
\end{figure}

(iii)
Among functions $c_i$, $i=1,2,3$ depending on $(x_*,y_*)$ the function $c_2(x_*,y_*)$ corresponds to the coefficient of $\alpha_1$ in the asymptotics and vanishes on a horizontal interface.

\vspace{1cm}
(iv) We observe the following values for Mumford-Shah energy $MS(x_*,y_*)$
and $c_2(x_*,y_*)$ in the minimum point $(-0.02,-0.02)$ and in the point $(0.000145,-0.01060)$
of intersection of the vertical interface  $\{SIF(x_*,y_*)=1\}$
and horizontal interface $\{ c_2(x_*,y_*)=0 \}$

\begin{table}[hptb]
\begin{center}
\begin{tabular}{|c|c|c| } \hline
   $ $   & $(-0.02,-0.02)$ & $(0.000145,-0.01060)$ \\ \hline
  $MS(\cdot,\cdot)$ & $3.2523$ & $3.2531$ \\ \hline
  $c_2(\cdot,\cdot)$ & $0.001689686245$ & $0$ \\ \hline
\end{tabular}
\end{center}
\end{table}

The small difference of the values and the deviation of the both points from the origin in the same direction, due to the perturbation of the boundary data with $\epsilon=0.01$, indicates that the conjecture stated in \cite{AM} that in the asymptotic expansions (\ref{asymptcrack}) and (\ref{asymptmt}) one should have $k>1$, could be the missing Euler-Lagrange condition.

\begin{figure}[htbp]\label{fig:c2}
\begin{minipage}[t]{2.3in}
\includegraphics[scale=0.38]{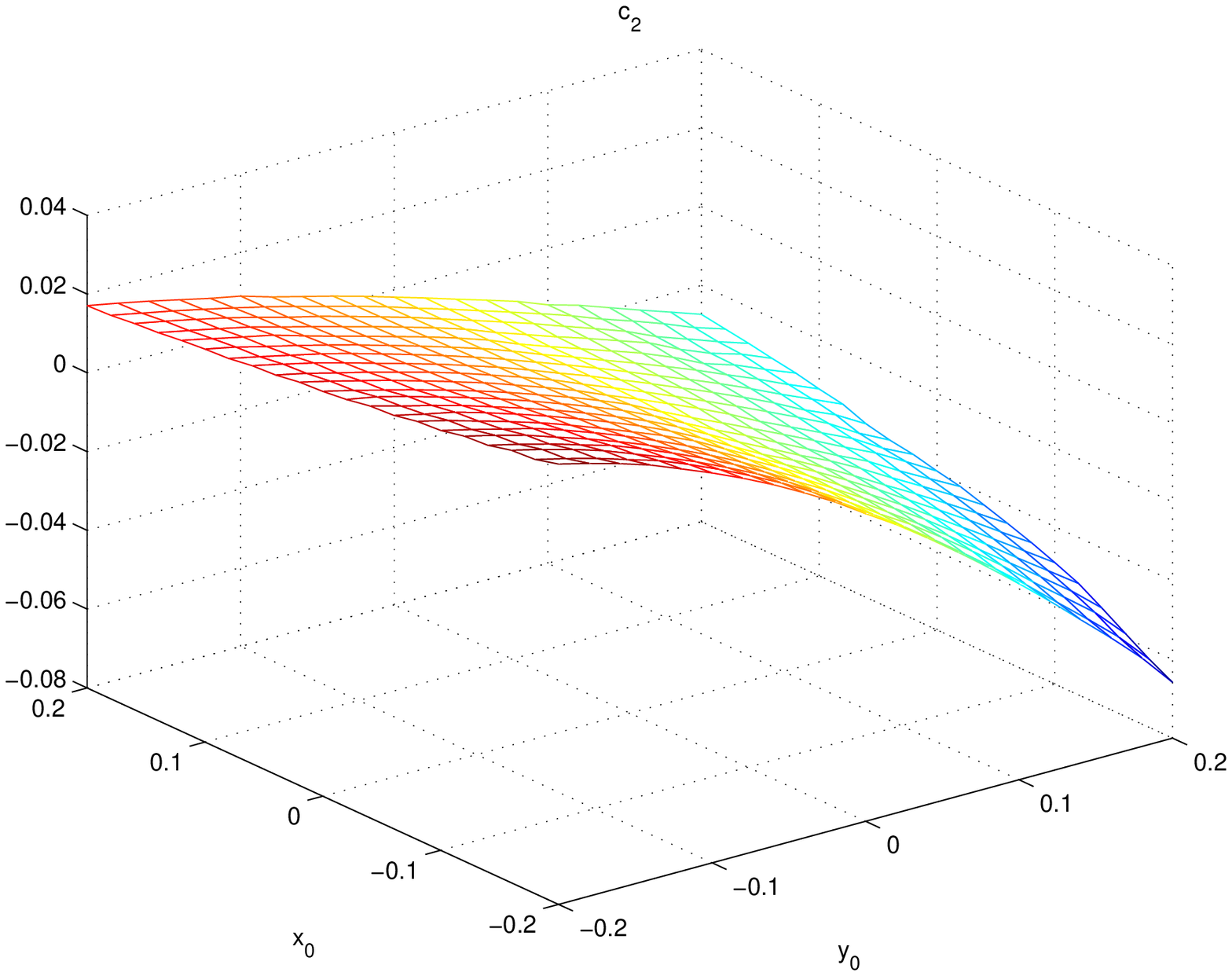}
\end{minipage}
\hfil
\begin{minipage}[t]{2.1in}
\includegraphics[scale=0.33]{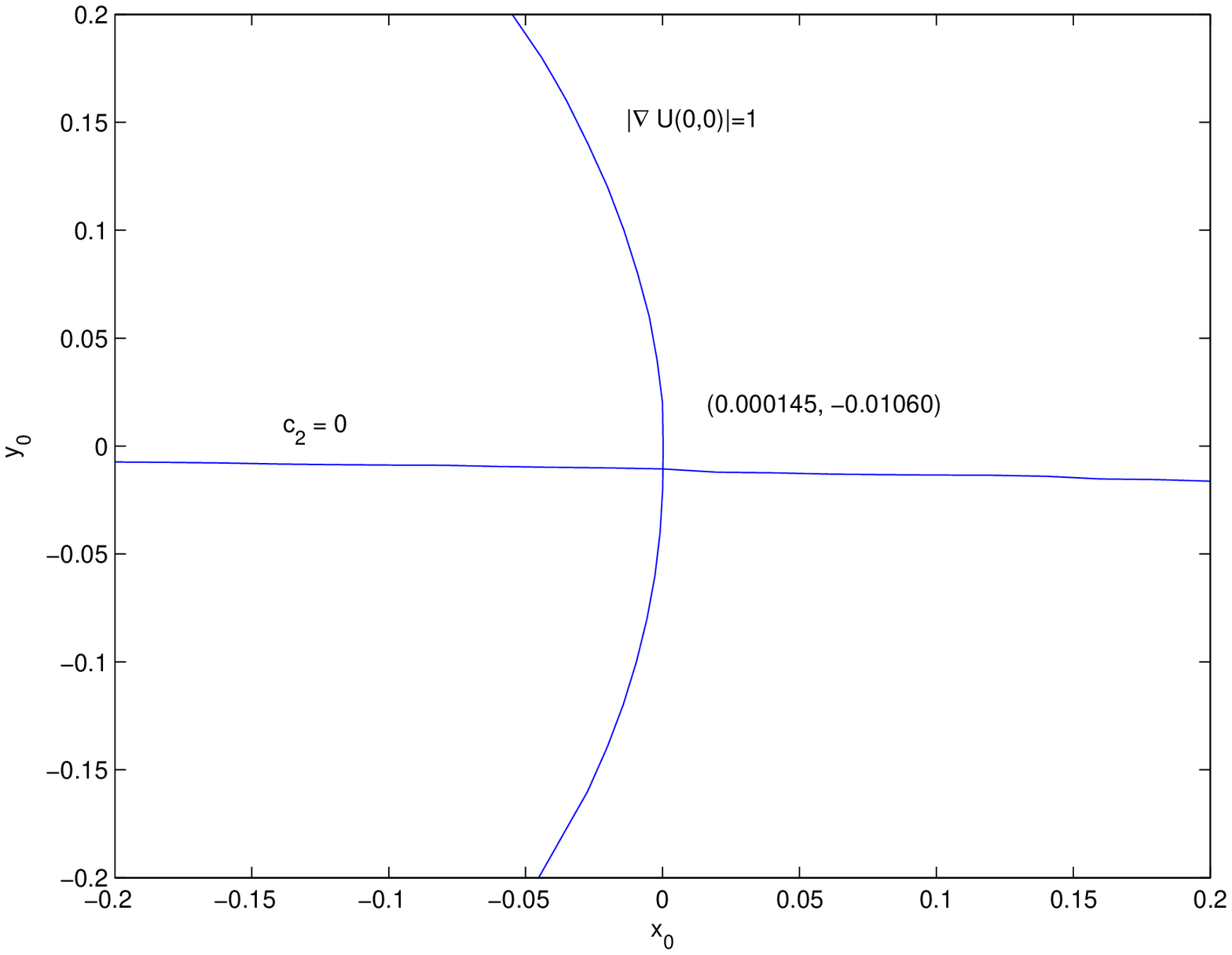}
\end{minipage}
\caption{Plot of the $c_2$ coefficient (left), and intersection of interfaces.} \label{fig:cross}
\end{figure}

\vspace{1cm}

(v)
We also notice that the value of the minimizer function $u$ on the crack-tip is changing linearly when it moves away from the minimizing point
(see Figure~\ref{fig: U00_UT00}, left).
This is an important observation, which might explain why all attempts to detect a missing first order condition by
traditional domain variation techniques have not yet led to any result.

\begin{figure}[htbp]\label{fig:c1}
\begin{minipage}[t]{2.3in}
\includegraphics[scale=0.35]{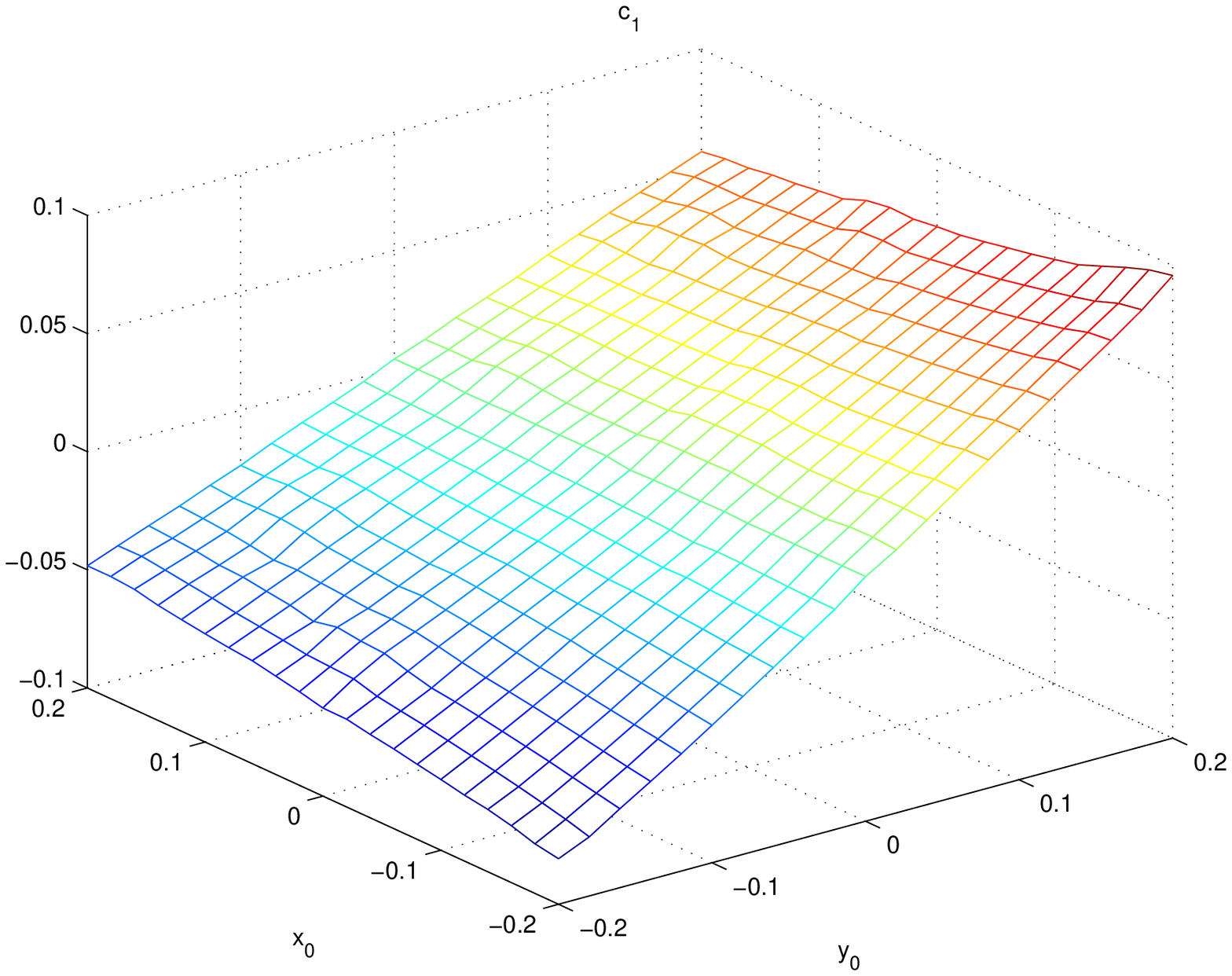}
\end{minipage}
\hfil
\begin{minipage}[t]{2.1in}
\includegraphics[scale=0.35]{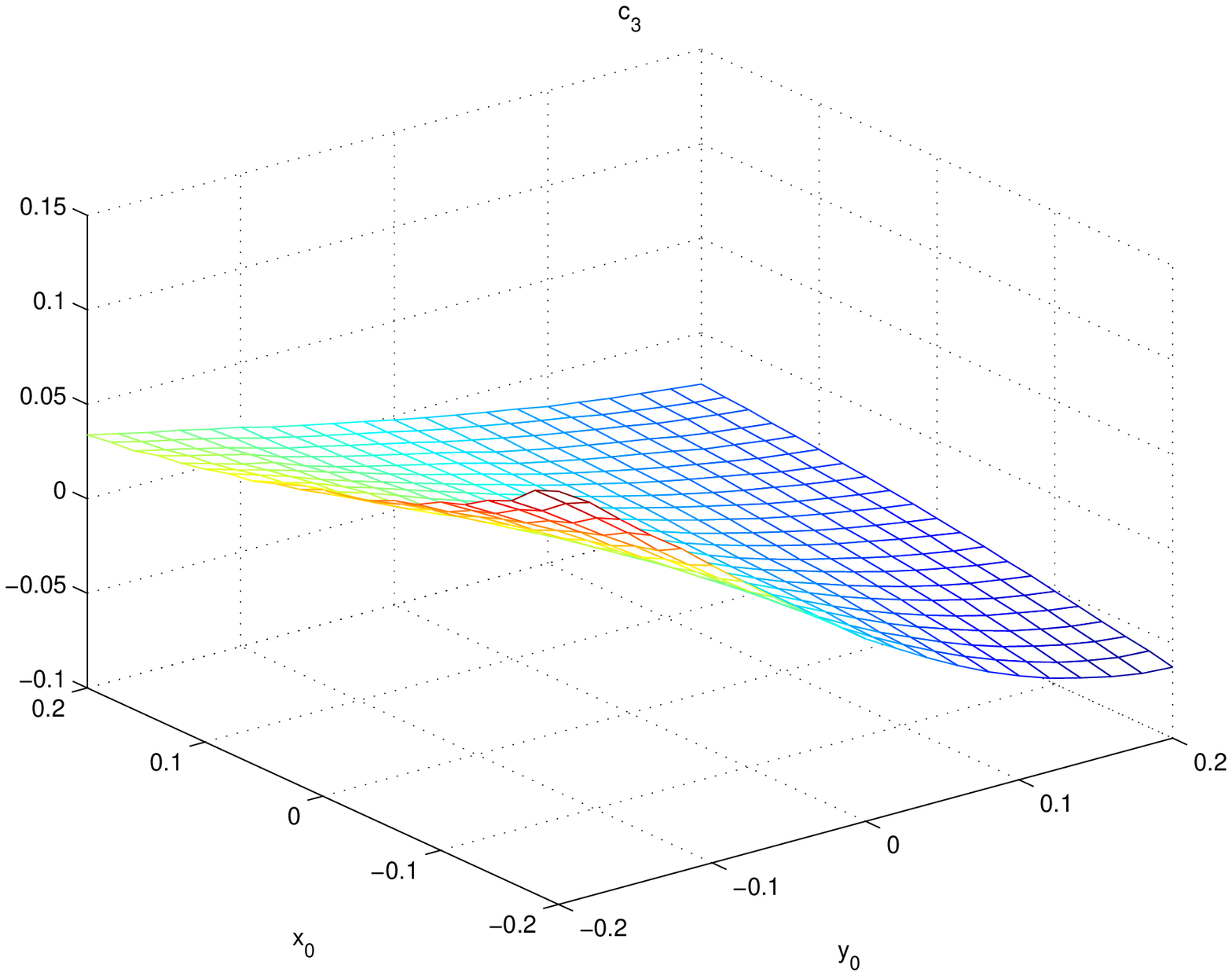}
\end{minipage}
\caption{Plot of the $c_1$ coefficient (left), and the $c_3$ coefficient (right).} \label{fig:c3}
\end{figure}

\subsection*{Acknowledgements}
The first author was partially supported by the US NSF grant DMS-1522768, and the NIH grant 5R01GM96195-2, CNSF grants 11371199, 11471166, and BK20141443. The second author was partially supported by the XJTLU's research development grant RDF-13-02-13.

The second author is grateful to
Stephan Luckhaus and Charles Elliott for inspiring discussions, as well as to Henrik Shahgholian and the Isaac Newton Institute in Cambridge for warm hospitality and
friendly atmosphere during the Research Program on Free Boundary Problems and Applications.


\end{document}